\documentclass[11pt]{amsart}
\usepackage{amsmath}
\usepackage{palatino}
  \usepackage{paralist}
  \usepackage{graphics} 
  \usepackage{epsfig} 
\usepackage{graphicx}  \usepackage{epstopdf}
\usepackage[pagewise]{lineno}

\newtheorem{theorem}{Theorem}[section]
\newtheorem{corollary}{Corollary}

\newtheorem{lemma}[theorem]{Lemma}
\newtheorem{proposition}{Proposition}

\theoremstyle{definition}

\newtheorem{remark}{Remark}

\newcommand{\R}{\mathbb{R}}
\newcommand{\C}{\mathbb{C}}

\newcommand{\be}{\begin{equation}}
\newcommand{\ee}{\end{equation}}

\title[ Klein-Gordon Schr\"{o}dinger system]
      {Global Attractor for Weakly damped nonlinear parametric Klein-Gordon Schr\"{o}dinger system}

\author[Olivier Goubet and Marilena Poulou]{}


 \email{olivier.goubet@univ-lille.fr}
 \email{  mpoulou@uniwa.gr}

\begin{document}
\maketitle

\centerline{\scshape Olivier Goubet}
\medskip
{\footnotesize
\centerline{Laboratoire Paul Painlevé CNRS UMR 8524,}

\centerline{Université de Lille, INRIA,}
\centerline{59 655 Villeneuve d'Ascq cedex.}

\medskip
\centerline{and }

} 

\medskip

\centerline{\scshape  Marilena Poulou}
{\footnotesize
\centerline{Department of Mechanical Engineering,}

\centerline{University of West Attica,,}
\centerline{Egaleo, 12241, Athens.}

\bigskip

 \centerline{\it Dedicated friendly to V. Pata, a fine mathematician, for his 60th birthday}

\begin{abstract}
We consider a damped forced nonlinear Klein-Gordon Schr\"{o}dinger (KGS) system
known as the parametric KGS system with Yukawa coupling.
We prove the regularity of the global attractor for this infinite-dimensional dynamical
system, i.e. that the global attractor is included into a smaller Sobolev space
than the energy space. We prove that this global attractor has finite fractal dimension,
without any assumption on the decay of the forcing terms.
\end{abstract}

\section{Introduction}

We consider a Klein-Gordon Schr\"{o}dinger system (KGS) with Yukawa coupling.
This system describes the nonlinear interaction between high frequency
electron waves and low frequency ion plasma waves in a homogeneous magnetic field,
adapted to model the upper hybrid plasma
heating scheme (see \cite{KSX} and the references therein).
This system couples a Schr\"odinger equation with
a wave equations and reads
\begin{eqnarray}
 \label{e1.1} \psi_t +i\Delta \psi  +\alpha\psi &= & i\phi \psi +f, \\
 \label{e1.2}
 \phi_{tt}-\Delta \phi +\phi+\gamma\phi_t  &= &  -Re \psi_{x} +g.
\end{eqnarray}

\noindent Here the unknowns $\psi$ and $\phi$ map respectively
$\R_t\times \R_x$ into respectively $\C$ and $\R$. The positive parameters $\alpha$
and $\gamma$ are the damping parameters that are supposed to be small, i.e. $\alpha, \gamma<<1$. In the right hand side $f, g$ are external forces
which do not depend on $t$ and that belong to $L^2(\R)$.

This weakly damped forced KGS system partakes of infinite-dimensional
dynamical systems  \cite{Hale}, \cite{lady}, \cite{MZ},\cite{Raugel}, \cite{Temam}.
This system is dissipative. Any trajectory enters a bounded subset
of the energy space $E=H^1(\R)\times H^1(\R)\times L^2(\R)$ after a finite transient time.
For this system, the second author has proved the existence
of a global attractor $\mathcal{A}$ in dimension $1$ \cite{PS1} \cite{PS2}
in the energy space $E$; a global attractor is a compact subset of $E$ that is invariant
by the flow of solutions and that attracts all trajectories when $t\rightarrow +\infty$, uniformly on bounded sets of $E$.
Besides, there are other Klein-Gordon-Schr\"odinger systems with Yukawa coupling
that were studied in the last decades. The system with $-|\psi|^2$ in the right hand side
of \eqref{e1.2} instead of -Re$\psi_x$ where the existence of the attractor
is known \cite{AGH}, \cite{GuoLi}, \cite{LuWang}, \cite{WangLange};
for this equation in the conservative case $\alpha=\gamma=0$ and $f=g=0$ the following Hamiltonian/Energy is conserved
along the flow of solutions
$$ \int (\phi_t)^2+\phi^2+|\nabla \phi|^2+2\phi|\psi|^2+|\nabla \psi|^2)dx.$$
The Zakharov system with $\Delta |\psi|^2$ instead of -Re$\psi_x$ was studied
in \cite{Flahaut}, \cite{GoubetMoise}; for this system in the conservative case the Hamiltonian reads, for $\Lambda=(-\Delta)^{-\frac12}$,
$$ \int (\Lambda \phi_t)^2+(\Lambda \phi)^2+\phi^2+2\phi|\psi|^2+|\nabla \psi|^2)dx.$$
In other words, these two systems are weakly damped perturbations
of Hamiltonian systems. The system under consideration in the present article is not.
Moreover, the second equation \eqref{e1.2} is linear but involves a derivative.

Our first result concerns the regularity of the global attractor.
The initial value problem related to \eqref{e1.1}-\eqref{e1.2} is well-posed
in the energy space $E$. Since we can go backward in time for the initial value problem
the flow map $S(t): (\psi(0), \phi(0),\phi_t(0))\mapsto (\psi(t), \phi(t),\phi_t(t))$
has no smoothing effect. Some regularization occurs when $t=+\infty$; this asymptotical smoothing effect reads as follows

\begin{theorem}\label{th1}
The compact global attractor $\mathcal{A}$ is a compact subset of $H^2(\R )\times H^2(\R )\times H^2(\R).$
\end{theorem}

This result is sharp since the global attractor contains the stationary solutions of the equation, i.e.
solution $(\psi^*,\phi^*,0)$ that solves
\begin{eqnarray}
 \label{e.1} i\Delta \psi^*  +\alpha\psi^* &= & i\phi^* \psi^* +f, \\
 \label{e.2}
 -\Delta \phi^* +\phi^* &= &  -Re \psi_{x}^* +g.
\end{eqnarray}

\noindent These stationary solutions belong to $H^2(\R )\times H^2(\R )\times \{0\}.$
This asymptotical smoothing effect was observed in \cite{haraux}, \cite{ghidagliatemam} for wave equations
and in \cite{G1} for Schr\"odinger equations.

Our second result asserts than the global attractor has finite fractal (or box) dimension.
Besides, when considering the issue of the fractal dimension of the global attractor,
it is usually assumed that the external forcing term possesses some extra decay property;
see for instance \cite{akroune} for the corresponding result for the one dimensional
Schr\"odinger equation on the whole line, with the assumption $\int_\R x^2|f(x)|^2dx<\infty$.
Our second result states

\begin{theorem}\label{th2}
The fractal dimension of $\mathcal{A}$ in $H^1(\R)\times H^1(\R)\times L^2(\R)$ is finite.
\end{theorem}

\begin{remark}
We point out that we do not need any decay assumption on $f,g$.
\end{remark}

The remaining of the article is organized as follows. In
Section 2 we prove Theorem \ref{th1}. The proof outlines as follows.
Dealing with the first equation we prove that the $\psi$ part of a trajectory
in the global attractor is more regular than expected; then we address the $(\phi,\phi_t)$
part and we bootstrap the argument. In Section 3 we prove Theorem \ref{th2} combining the method
of contraction of infinitesimal volumes (see \cite{Temam} and the references therein)
and a new approach introduced by B. Alouini \cite{Alouini} using restriction norms .
In Section 4, we give an example of Parametric KGS system whose global attractor is not trivial,
i.e. is not reduced to a single point.

\smallskip
To complete this introduction, let us introduce some notations. The Hilbert space $L^2(\R)$ is equipped with the scalar product
 $$ {\rm Re}\int_{\mathbb{R}}u(x) \overline{v}(x)dx.$$
We also use standard Sobolev and Banach spaces as $H^1(\R)$ or $L^\infty(\R)$. If $p$ is a nonnegative function in $\R$
we define the space $$L^2(p(x)dx)=\{ f \in L^1_{loc(\R)}; \int_\R p(x) |f(x)|^2dx <+\infty\}.$$
The set $L^2(-R,R)$ is the set of square integrable functions defined on the interval $[-R,R]$.
\noindent By the letter $c$ we denote a numerical constant that may change from one line to one
another. The letter $K$ will denote some generic constant that
depends on the data $\alpha, \gamma, f, g$.
For instance the global attractor $\mathcal{A}$ is included in a ball of radius $K$
in the energy space.

\section{Proof of Theorem \ref{th1}}

Let us recall that a {\it complete} trajectory $(\psi(t),\phi(t),\phi_t(t))$
is a solution of \eqref{e1.1}-\eqref{e1.2} that is defined for all $t$ in $\R$
and that is bounded in the energy space $E$. Let us recall then
that the global attractor $\mathcal{A}$ is the set of complete trajectories
for the Klein-Gordon-Schr\"odinger system \cite{Temam}.

\subsection{Extra regularity on $\psi$}
\label{firststep}

We prove here that the first component $\psi(t)$ of
a  complete bounded trajectory $(\psi(t),\phi(t),\phi_t(t))$ in the energy space
is in fact included in a bounded set of $H^2(\R)$.

Let $m$ be in $\mathbb{N}$. We split $\psi(t)=u^m(t)+v^m(t)$ where $u^m$ is solution for $t\geq -m$
to the equation \eqref{mp1} below with non standard initial condition at $t=-m$
that reads $u^m(-m)=u_m$ where

\begin{eqnarray}
\label{mp1} u^m_t +i\Delta u^m  +\alpha u^m &= & i\phi u^m+f, \\
\label{mp1.0}\alpha u_m+i \Delta u_m &= &i\phi(-m)u_m+f.
\end{eqnarray}

\noindent For a pair $(\phi(-m),f)$ there exists a unique
$u_m$ solution of the elliptic equation \eqref{mp1.0} thanks to Lax-Milgram theorem.
Then the remaining term $v^m=\psi-u^m$ solves

\begin{equation}\label{mp2}
 v^m_t +i\Delta v^m  +\alpha v^m =  i\phi v^m ,
\end{equation}

\noindent with initial data $v^m(-m)=v_m=\psi(-m)-u_m$ at $t=-m$.
We now establish that $v^m$ vanishes when $m$ goes to $+\infty$.
Considering the scalar product of \eqref{mp2} with $v$ and integrating in $t$ between $-m$ and $0$ leads to

 \begin{equation}\label{mp3}
 ||v^m(0)||^2_{L^2(\R)}\leq ||v^m(-m)||^2_{L^2(\R)}\exp(-2\alpha m)\leq K \exp(-2\alpha m).
 \end{equation}

\noindent Considering the scalar product of \eqref{mp2} with $-\Delta v^m$ leads to,
since $\phi_x$ remains in a bounded set of $L^2(\R)$, and using the sharp Gagliardo-Nirenberg inequality
(see \cite{HolmerLiu} for a proof)
$$ ||u||_{L^\infty(\R)}^2\leq ||u||_{L^2(\R)}||u_x||_{L^2(\R)},$$
\noindent yields

 \begin{equation}\label{mp3.5}
 \frac12 \frac{d}{dt} ||v^m_x||^2_{L^2(\R)}+\alpha ||v^m_x||^2_{L^2(\R)}=-{\rm Im}\int_\R \phi_x v^m \overline{v}^m_x\leq K ||v^m||^\frac12_{L^2(\R)} ||v^m_x||^\frac32_{L^2(\R)}.
 \end{equation}

 \noindent Then thanks to Young inequality and to Gronwall lemma

 \begin{equation}\label{mp3.6}
 ||v^m_x(t)||^2_{L^2(\R)}\leq e^{-\alpha (t+m)}||v^m_x(-m)||^2_{L^2(\R)}+K\int_{-m}^t e^{\alpha(s-t)} ||v^m(s)||^2_{L^2(\R)} ds.
 \end{equation}

\noindent Therefore  $v^m_x$ and then $u^m_x$ remain in a bounded set of $H^1(\R)$,
and $v^m(0)$ converges to $0$ when $m$ diverges to $+\infty$.

 We now prove that $u^m$ belongs to a more regular space than $H^1(\R)$.
 Differentiating \eqref{mp1} with respect to $t$ and setting $w^m=u^m_t$ leads to, thanks to the initial condition \eqref{mp1.0}

 \begin{eqnarray}
\label{mp4} w^m_t +i\Delta w^m  +\alpha w^m &= & i\phi w^m+i \phi_t u^m, \\
\label{mp4.0}w^m(-m)&= &0.
\end{eqnarray}

Considering the scalar product of \eqref{mp4} by $w^m$ leads to

 \begin{equation}\label{mp5}
 \frac12 \frac{d}{dt} ||w^m||^2_{L^2(\R)}+\alpha ||w^m||^2_{L^2(\R)}=-{\rm Im}\int_\R \phi_t u^m \overline{w}^m.
 \end{equation}

\noindent Since the right hand side of \eqref{mp5} is bounded by

$$||w^m||_{L^2(\R)} ||\phi_t||_{L^2(\R)} ||u^m||_{L^\infty(\R)}\leq K ||w^m||_{L^2(\R)},$$

\noindent applying Gronwall lemma leads to $||w^m(0)||_{L^2(\R)}\leq K \alpha^{-1}.$

\noindent Thanks to \eqref{mp1} we have $||u^m(0)||_{H^2(\R)}\leq K$.
Then $\psi(0)$ splits $\psi(0)=u^m(0)+v^m(0)$, where the sequence $m\mapsto u^m(0)$ is bounded in $H^2(\R)$,
and where $v^m(0)$ converges towards $0$ in $L^2$ when $m$ diverges to $+\infty$.
Then we have a uniform bound in $H^2(\R)$ for $\psi(0)$.

Using the same argument replacing $t=0$
by any $t$ yields that the first component $\psi(t)$ of a complete trajectory remains in a bounded set of $H^2(\R)$.

\subsection{Extra regularity on $\phi,\phi_t$}

The $H^2(\R)\times H^1(\R)$ bound on $(\phi,\phi_t)$
is a consequence of the following lemma, with
$\sigma=1$, $u=\phi-\phi^*$ and $h={\rm Re}(\psi^*-\psi_x)$.

 \begin{lemma}\label{lemma} Consider $\sigma$ in $\R$.
Consider a function $h(t,x)$ that is continuous and bounded from $\R_t$
into $H^{\sigma}(\R)$. Consider $(u,u_t)$ the solution to
\be\label{wave}
u_{tt}+\gamma u_t-\Delta u+u=h,
\ee
\noindent that is continuous and bounded from $\R_t$ into $H^\sigma(\R)\times H^{\sigma-1}(\R)$.
Then $(u,u_t)$ is a continuous and bounded function from $\R_t$ into $H^{\sigma+1}(\R)\times H^{\sigma}(\R)$.
\end{lemma}

\begin{proof}

\noindent Following \cite{Bourgain}, we introduce the new variable $v=u+i(Id-\Delta)^{-\frac12}u_t$
that is solution to the Schr\"odinger like equation

\begin{equation}\label{nls}
v_t+i\gamma {\rm Im} v+i(Id-\Delta)^\frac12 v=i (Id-\Delta )^{-\frac12}h.
\end{equation}

\noindent Consider the family of linear operators
defined by $U_\gamma(t)v_0=v(t)$ where $v$ solves
\eqref{nls} with $h=0$, and with $v(0)=v_0$.
Using Fourier transform in $x$, and setting $$\hat{v}(t)(\xi)=Y(t,\xi),$$
\noindent we have that

\begin{equation}\label{can1}
\frac{dY}{dt} +\begin{pmatrix} & 0  &-\sqrt{1+\xi^2} \\ & \sqrt{1+\xi^2} & \gamma\end{pmatrix} Y= \frac{dY}{dt}  + L(\xi) Y=0.
\end{equation}

\noindent Actually

\begin{equation}\label{can2}
e^{-tL(\xi)}=P^{-1}(\xi) \begin{pmatrix}  & e^{-t\lambda_-(\xi)} & 0\\  &0 & e^{-t\lambda_+(\xi)} \end{pmatrix} P(\xi),
\end{equation}

\noindent where $||P(\xi)||\;|| P^{-1}(\xi)||\leq K$ and $2\lambda_\pm(\xi)=\gamma \pm i\sqrt{4(1+\xi)^2-\gamma^2}$.
Since $\gamma^2<4+4\xi^2$ then the two eigenvalues are complex conjugate and
$||U_\gamma(t)||\leq K\exp(-\frac \gamma 2 t)$.
Therefore the only bounded solution of the equation
in the statement of Lemma \ref{lemma} is

\begin{equation}\label{rep}
v(t)=i\int_{-\infty}^t U_{\gamma}(t-s) (Id-\Delta)^{-\frac12} g(s)ds.
\end{equation}

\noindent The proof of the Lemma is then a mere consequence
of this representation formula. For instance

\begin{equation}\label{final}
||v(t)||_{H^{\sigma+1}(\R)}\leq \left(\int_{-\infty}^t || U_{\gamma}(t-s)||ds\right)\sup_s ||h||_{H^{\sigma}(\R)}\leq  \frac 2 \gamma \sup_s ||h||_{H^{\sigma}(\R)}.
\end{equation}

\end{proof}

\begin{remark}
In fact the representation formula \eqref{rep} shows that the term $i\phi \psi$
acts as a nonlinear infinite memory term for the dynamics in the attractor.
\end{remark}

\subsection{Extra regularity on $\psi_t$}
\label{thirdstep}

Consider a complete trajectory in the global attractor.
The equation \eqref{e1.1} and the two previous subsections yield that $\psi_t$
remains in a bounded set of $L^2(\R)$. In fact this result can be improved

\begin{proposition}\label{inexpected}
Consider a complete trajectory in the global attractor. Then $\psi_t$ remains bounded in $H^1(\R)$.
\end{proposition}

\begin{proof}
We use the splitting $\psi=u^m+v^m$ introduced in Subsection \ref{firststep}.
The part $v^m$ is the vanishing part when $m$ diverges towards $+\infty$;
 the derivative $v^m_t$ vanishes also thanks to \eqref{mp2}.
We focus on $w^m=u^m_t$. Considering the scalar product of \eqref{mp4} with $-i(w^m_t+\alpha w^m)$ leads
to, for $t\geq -m$

\begin{equation}\label{in1}
\frac 1 2 \frac{d}{dt} J(w^m)+\alpha J(w^m)=\frac 32 \int_\R \phi_{t} |w^m|^2 +{\rm Re} \int_\R \phi_{tt} u^m \overline{ w}^m+
\alpha {\rm Re} \int_\R \phi_t u^m \overline{w}^m ,
\end{equation}

\noindent where $J(w^m)=||w^m_x||^2_{L^2(\R)}+\int_\R \phi |w^m|^2 +2{\rm Re}\int_\R \phi_t u^m \overline{w}^m.$
Due to the previous estimates on the global attractor

\begin{equation}\begin{split}\label{in2}
J(w^m) &\geq  ||w^m_x||^2_{L^2(\R)}- ||\phi||_{L^\infty} ||w^m||^2_{L^2(\R)}- 2||\phi_t||_{L^2(\R)}||u^m||_{L^\infty(\R)} ||w^m||_{L^2(\R)}\\  &\geq  ||w^m_x||^2_{L^2(\R)}-  K.
\end{split}\end{equation}

\noindent We now handle the right hand side of \eqref{in1}. We have

\begin{equation}\label{in3}
\frac 12 \int_\R \phi_{t} |w^m|^2 \leq ||\phi_{t}||_{L^\infty(\R)}||w^m||^2_{L^2(\R)}\leq K.
\end{equation}

\noindent Analogously

\begin{equation}\label{in5}
|\int_\R \phi_{tt} u^m \overline{ w}^m|\leq ||\phi_{tt}||_{L^2(\R)}||u^m||_{L^\infty(\R)}||w^m||_{L^2(\R)}\leq K,
\end{equation}

\noindent and

\begin{equation}\label{in7}
|\int_\R \phi_t u^m \overline{w}^m| \leq ||\phi_{t}||_{L^2(\R)}||u^m||_{L^\infty(\R)}||w^m||_{L^2(\R)}\leq K.
\end{equation}

\noindent Gathering these inequalities the Gronwall Lemma completes the proof of the Proposition.

\end{proof}

\begin{corollary}\label{extra2}
Consider a complete trajectory in the global attractor. Then $\phi_t$ remains bounded in $H^2(\R)$.
\end{corollary}

\begin{proof}
Differentiate \eqref{e1.2} with respect to $t$ and apply Lemma \ref{lemma} with $\sigma=1$ to
the equation

\begin{equation}\label{extra3}
(\phi_t)_{tt}+\gamma (\phi_t)_t+\phi_t-\Delta \phi_t=-{\rm Re}(\psi_{tx}).
\end{equation}

\end{proof}

\subsection{Compactness of the global attractor in $H^2(\R)\times H^2(\R) \times H^2(\R)$}

We first recall the following lemma

\begin{lemma}\label{compact}
Consider a closed set $\mathcal{E}$ in Banach space $E$ such that for any $\varepsilon>0$
there exists a compact set $K_\varepsilon$ of $E$ such that $\mathcal{E}\subset K_\varepsilon+B_E(0,\varepsilon)$,
where $B_E(0,\varepsilon)$ is the closed ball in $E$ centered at $0$ and of radius
$\varepsilon$. Then $\mathcal{E}$ is a compact subset of $E$.
\end{lemma}

\begin{proof}
For any $\varepsilon>0$, due to the assumptions, the set $\mathcal{E}$ can be covered
by a finite number of balls of radii $\varepsilon$; in other words,
$\mathcal{E}$ is totally bounded. Since it is a complete metric space, then it is compact.
\end{proof}

We now state a compactness result for the set
of the first component $\psi$ of points in the attractor.
More precisely this set is the set of the first component
of a complete trajectory evaluated at $t=0$.

\begin{lemma}\label{pan1}
The set $\{ \psi(0) ; (\psi(t),\phi(t), \phi_t(t)) \in \mathcal{A}\}$
is a compact set of $L^2(\R)$.
\end{lemma}

\begin{proof}

We use here a splitting of the solution as in the proof of the existence
of a global attractor in \cite{Laurencot}.
For $\varepsilon>0$ consider a compactly supported function $f_\varepsilon$ in $H^1(\R)$ such that
$||f-f_\varepsilon||_{L^2(\R)}\leq \varepsilon$. Consider now the splitting $\psi=u^{m,\varepsilon}+v^{m,\varepsilon}$
where for $t\geq -m$,

\begin{eqnarray}
\label{mf1} u^{m,\varepsilon}_t +i\Delta u^{m,\varepsilon}  +\alpha u^{m,\varepsilon} &= & i\phi u^{m,\varepsilon}+f_\varepsilon, \\
\label{mf1.0}u^{m,\varepsilon}(-m)&= &0.
\end{eqnarray}

\noindent and

\begin{eqnarray}\label{mf2}
 v^{m,\varepsilon}_t +i\Delta v^{m,\varepsilon}  +\alpha v^{m,\varepsilon} &= & i\phi v^{m,\varepsilon}+f-f_\varepsilon, \\
\label{mf2.0} v^{m,\varepsilon}(-m)&= &\psi(-m).
\end{eqnarray}

\noindent Considering the scalar product of \eqref{mf2} with $ v^{m,\varepsilon} $ yields

\begin{equation}\label{mf3}
\frac12 \frac{d}{dt} ||v^{m,\varepsilon}||^2_{L^2(\R)}+\alpha  ||v^{m,\varepsilon}||^2_{L^2(\R)}\leq ||f-f_\varepsilon||_{L^2(\R)} ||v^{m,\varepsilon}||_{L^2(\R)},
\end{equation}

\noindent and then by Gronwall's Lemma

\begin{equation}\label{mf4}
||v^{m,\varepsilon}(0)||_{L^2(\R)}\leq \exp(-\alpha m) ||\psi(-m)||_{L^2(\R)}+\frac \varepsilon \alpha.
\end{equation}

\noindent Similarly to the computations in Subsection \ref{firststep} we can prove that $v^{m,\varepsilon}$
remains in a bounded set of $H^1(\R)$. Then $u^{m,\varepsilon}$ remains in a bounded set of $H^1(\R)$ as well.

\noindent Considering the scalar product of \eqref{mf1} with $x^2 u^{m,\varepsilon}$ yields

\begin{equation}\begin{split}\label{mf3}
\frac12 \frac{d}{dt} ||xu^{m,\varepsilon}||^2_{L^2(\R)}+\alpha  ||xu^{m,\varepsilon}||^2_{L^2(\R)}&\leq  |{\rm Im}\int_\R 2x u^{m,\varepsilon}_x \overline{u}|
+|{\rm Re}\int_\R x^2f_\varepsilon \overline{u}^{m,\varepsilon}dx|\\  &\leq  \frac \alpha 2 ||xu^{m,\varepsilon}||^2_{L^2(\R)}+ K+ ||xf_\varepsilon||_{L^2(\R)}^2.
\end{split}\end{equation}

\noindent Thanks to Gronwall's lemma we have that $xu^{m,\varepsilon}(0)$ remains in a bounded set of $L^2(\R)$,
the bound depends on $\varepsilon$ but not on $m$. Therefore letting $m$ diverge to $+\infty$
we have that $\psi(0)$ is captive in $K_\varepsilon+B_E(0,\frac \varepsilon \alpha)$ where $K_\varepsilon$
is closed and bounded in $H^1(\R) \cap L^2(x^2dx)$ and then compact in $L^2(\R)$.
\end{proof}

By interpolation with the $H^2(\R)$ bound  we have

\begin{corollary}
Set $\sigma<2$. The set $\psi=\psi(0)$ such that $(\psi(t),\phi(t), \phi_t(t))$ belongs to $\mathcal{A}$
is a compact set of $H^\sigma(\R)$.
\end{corollary}

We now propagate the compactness to the $\phi,\phi_t$ part of
the attractor.

\begin{lemma}\label{pan2}
For $\sigma<2$, the set $(\phi,\phi_t)=(\phi(0),\phi_t(0))$ such that $(\psi(t),\phi(t), \phi_t(t))$ belongs to $\mathcal{A}$
is a compact set of $H^\sigma(\R)\times H^{\sigma}(\R)$.
\end{lemma}

\begin{proof}
We use the representation formula \eqref{rep} that reads

\begin{equation}\label{rep2}
\phi(0)+i(Id-\Delta)^{-\frac12}\phi_t(0)=i\int_{-\infty}^0 U_{\gamma}(-s) (Id-\Delta)^{-\frac12}(-{\rm Re}\psi_x(s)+g)ds.
\end{equation}

Consider $\sigma\geq 1$. The map $\psi(0) \mapsto \psi(s)$ is continuous from $H^{\sigma}(\R)$
to $C([-T, 0]; H^\sigma(\R))$. This is a byproduct of solving the initial value problem.
Fix now $T$ large enough such that, thanks to \eqref{final}

$$ ||\int_{-\infty}^{-T} U_{\gamma}(-s) (Id-\Delta)^{-\frac12}(-{\rm Re}\psi_x(s)+g)ds||_{H^\sigma(\R)}\leq \varepsilon.$$

\noindent This part is then included in a small ball.
\noindent Besides, since the map $$F: \psi(0)\mapsto \int_{-T}^0 U_{\gamma}(-s) (Id-\Delta)^{-\frac12}(-{\rm Re}\psi_x(s)+g)ds$$
\noindent is continuous, then $F$ maps a compact set into a compact set.
Thanks to Lemma \ref{compact} the pair $(\phi(0),\phi_t(0))$ remains in compact set of $H^\sigma(\R)\times H^{\sigma-1}(\R)$.
For the $H^\sigma(\R)$ compactness we interpolate with the $H^2(\R)$ bound for $\phi_t$.
\end{proof}

To reach the compactness of $\psi$ for $\sigma=2$, we prove that $w=\psi_t$
is in a compact set of $L^2(\R)$. Differentiating \eqref{e1.1} with respect to
$t$ and considering the scalar product of the resulting equation with $w$ yields

\begin{equation}\label{Ball}
||w(t)||^2_{L^2(\R)}=||w(0)||^2_{L^2(\R)}e^{-2\alpha t}-\int_0^t e^{2\alpha(s- t)}(\int_\R {\rm Im} \phi_t(s)\psi(s)\overline{w}(s)dx) ds.
\end{equation}

The right hand side is weakly continuous in $L^2(\R)$ since both $\phi_t$ and $\psi$
are compact in $H^1(\R)$. Then  the framework of the J. Ball's argument \cite{Ball} applies
\cite{mrw}, \cite{Wang} and we have the compactness of the $\psi$ part of the attractor in $H^2(\R)$.
For the $(\phi,\phi_t)$ part we propagate the compactness using Lemma \ref{pan2}.
This completes the proof of Theorem 1.1.

\section{The fractal dimension of the global attractor is finite}

In this section we prove that the global attractor has finite
fractal dimension.
 We follow here the Douay-Oesterl\'e \cite{DO} method
as described in \cite{Ghidaglia}, \cite{Temam} of contraction
of infinitesimal volumes. For alternate proofs that do not use the differentiability of the semigroup
we refer to \cite{zelik}, \cite{chueshov}.

\subsection{A general result}

Let $X$ be a separable Hilbert space.
Set respectively $(u,v)_X$ and $||u||_X$ for respectively the scalar product and norm in $X$.
For a collection $(u_1,...u_m)$ of $m$ vectors in $X$
set $G(u_1,....u_m)$ for the Gram determinant of this collection of vectors,
i.e. the determinant of the matrix $g$ whose entries are $g_{kl}=(u_k,u_l)_X$.
This quantity is the volume of the parallelepiped whose vertices are $(u_1,...u_m)$.
We then recall from \cite{Temam}, chapter $5$

\begin{theorem}\label{chulas}
Consider $(h_1,....h_m)$ an orthonormal collection of $m$ vectors in $X$.
Let $\mathcal{A}$ be a compact set in $X$ that is invariant by a semigroup $S(t)$ acting on $X$.
For $a$ in $\mathcal{A}$, consider $S(t)a$ the complete trajectory emanating from $a$
at $t=0$ and $u_j(t)=DS(t)(a)h_j$ the derivative of the semigroup emanating from $h_j$
at $t=0$. If the $m-$dimensional volume $G(u_1(t),....u_m(t))$ decays to $0$ when
$t$ goes to $+\infty$ then $\dim \mathcal{A}\leq m$.
\end{theorem}

We use the following corollary to conclude.
Let $Y$ be another space such that the embdding $X \subset Y$ is compact.
Set $(u,v)_Y$ and $||u||_Y$ for the scalar product and norm in $Y$.

\begin{corollary}
Let $\delta$ be a positive number.
If we have the following contraction estimate, up to a compact term
\begin{equation}\label{boa}
\frac{d}{dt} ||u||^2_X +\delta ||u||^2_X \leq K ||u||^2_Y,
\end{equation}
\noindent then the dimension of $\mathcal{A}$ is finite.
\end{corollary}

\begin{proof}
Since $X\subset Y$ is compact there exists a self-ajoint compact operator $T$
such that $(u,v)_Y=(u,Ty)_X$ for any vectors $u,v$ in $X$.
Moreover there exists a countable sequence of eigenvalues for $T$
that satisfies $\lambda_j\geq \lambda_{j+1}\rightarrow 0$.
Introduce $G(t)=G(u_1(t),..u_m(t))$. Then we infer from \eqref{boa} by the Courant min-max formula
(see \cite{Ghidaglia}, \cite{Temam} and the references therein)
\begin{equation}
\frac{d G}{dt} + 2\delta m G \leq K \left(\sum_{j=1}^m \lambda_j\right) G.
\end{equation}

\noindent Since $\frac 1 m \sum_{j=1}^m \lambda_j $ converges to $0$ when $m$
diverges to $+\infty$, then the volume $G(t)$ decays to $0$ when $t$ diverges to $+\infty$
if $m$ is such that $\frac 1 m \sum_{j=1}^m \lambda_j  < \frac{2\delta}{ K}$.
\end{proof}

\subsection{Proof of Theorem \ref{th2}}

We know (see \cite{PS2}) that the semi-group $S(t)$ is differentiable in the
energy space, for trajectories in the global attractor.

We follow now ideas introduced by B. Alouini (see \cite{Alouini} and the references therein).
The strategy reads as follows. Thanks to Rellich's Theorem, the embedding $H^1(\R) \subset L^2(-R,R)$ is compact.
Therefore the corresponding space $Y=L^2(-R,R)$  is suitable in estimates like \eqref{boa}.

The first step in Section \ref{firststepbeyond} below is then to obtain uniform bounds on the tails of the solution, that
are the semi-norms of a solution outside a given interval $[-R,R]$.
The second step in Section \ref{secondstepbeyond} consists in proving the key estimate \eqref{boa}
with $X$ the energy space and $Y=L^2(-2R,2R)\times \{0\} \times \{0\}$,
using that the map $(\psi,\phi,\phi_t) \mapsto (\psi,0,0)$ is compact.

\subsection{Estimates on the tails of a trajectory}
\label{firststepbeyond}

Consider a non negative smooth function $\theta$ such that $\theta(x)=1$ if $|x|\leq 1$ and $\theta(x)=0$ if $|x|>2$.
Set $\theta_R(x)=\theta(\frac x R)$. We define the seminorm representing the tail of $\psi$ as

\begin{equation}\label{seminorm}
||\psi||_{\theta}=(\int_\R (1-\theta_R(x))|\psi(x)|^2dx)^\frac12.
\end{equation}

\begin{lemma}\label{tail1}
For any $\varepsilon>0$ there exists a $R$ large enough depending on the bounds on the attractor such that
\begin{equation}\label{ta0}
||\psi||_\theta+ ||\psi_t||_\theta \leq \varepsilon.
\end{equation}
\end{lemma}

\begin{proof}

We begin with the $\psi$-estimate. Considering the scalar product of \eqref{e1.1} with $(1-\theta_R)\psi$ leads to

\begin{equation}\label{ta1}
\frac12 \frac{d}{dt} ||\psi||_{\theta}^2+\alpha ||\psi||_{\theta}^2= {\rm Im}\int_\R \psi_{xx}(1-\theta_R)\overline \psi +{\rm Re}\int_\R f(1-\theta_R)\overline \psi.
\end{equation}

\noindent The second term of the right hand side of \eqref{ta1} is bounded by above, thanks to Cauchy-Schwarz inequality, by
$||\psi||_{\theta} (\int_{|x|>R}|f|^2)^\frac12 $. Integrating by parts the first term in the right hand side of \eqref{ta1} we have

\begin{equation}\label{ta2}
| {\rm Im}\int_\R \psi_{xx}(1-\theta_R)\overline \psi|=|R^{-1}\int_\R \theta'(\frac x R) {\rm Im}( \psi_{x}(x)\overline \psi(x))dx|\leq R^{-1}||\theta'||_{L^\infty}||\psi||^2_{H^1(\R)}.
\end{equation}

\noindent Since $\psi$ is included in the global attractor and then remains bounded in $H^1(\R)$ we have

\begin{equation}\label{ta3}
\frac{d}{dt} ||\psi||_{\theta}^2+\alpha ||\psi||_{\theta}^2\leq \frac 1 \alpha \int_{|x|>R}|f|^2+ \frac K R.
\end{equation}

\noindent Integrating this between $-\infty$ and $0$ yields

\begin{equation}\label{ta3.5}
||\psi||_{\theta}^2\leq \frac{1}{\alpha^2}(\int_{|x|>R}|f|^2)+\frac{K}{\alpha R}.
\end{equation}

\noindent This completes the first part of the lemma, choosing $R$ large enough.
We now handle the $\psi_t$-estimate. Introduce now $w=\psi_t$ that is solution
to
\begin{equation}\label{orkyn}
w_t+\alpha w +iw_{xx}=i\phi w+i \phi_t \psi,
\end{equation}

\noindent that is defined for any $t$ in $\R$ and bounded in $H^1(\R)$.
Considering now the scalar product of \eqref{orkyn} with $(1-\theta_R) w$
leads to

\begin{equation}\label{ta4}
\frac12 \frac{d}{dt} ||w||_{\theta}^2+\alpha ||w||_{\theta}^2= {\rm Im}\int_\R w_{xx}(1-\theta_R)\overline w -{\rm Im}\int_\R \phi_t (1-\theta_R)\overline w \psi.
\end{equation}

\noindent The first term in the right hand side of \eqref{ta4} is bounded by above, proceeding as in \eqref{ta2}, by $R^{-1}||w||^2_{H^1}$.
The second term is bounded by $||\phi_t||_{L^\infty}||\psi||_\theta||w||_\theta$. Therefore

\begin{equation}\label{ta5}
\frac{d}{dt} ||w||_{\theta}^2+\alpha ||w||_{\theta}^2\leq R^{-1}||w||^2_{H^1(\R)}+\frac 1 \alpha ||\phi_t||_{L^\infty(\R)}^2||\psi||_\theta^2.
\end{equation}

\noindent Using Proposition \ref{inexpected}, the first part of the proof of this Lemma, and
integrating between $-\infty$ and $0$ completes the proof.
\end{proof}

We proceed similarly for the $(\phi,\phi_t)$.

\begin{lemma}\label{tail2}
For any $\varepsilon>0$ there exists a $R$ large enough depending on the bounds on the attractor such that
\begin{equation}\label{ta00}
||\phi||_\theta+||\phi_x||_\theta+ ||\phi_t||_\theta \leq \varepsilon.
\end{equation}
\end{lemma}

\begin{proof}

Let us recall that $0<\gamma<<1$.
Following an idea introduced in \cite{ghidagliatemam}, the equation \eqref{e1.2} reads also

\begin{eqnarray}
 \label{e1.21} \phi_t +\frac \gamma 2 \phi &= & y, \\
 \label{e1.22}
 y_{t}+ \frac \gamma 2 y -\Delta \phi +(1-\frac{\gamma^2}{4})\phi &= &  -Re \psi_{x} +g.
\end{eqnarray}

\noindent We consider the scalar product of the first equation by $(1-\theta_R)( -\Delta \phi +(1-\frac{\gamma^2}{4})\phi)$
and the scalar product of the second equation by $(1-\theta_R)y$. Summing the resulting equations and
setting $I= ||y||_{\theta}^2+||\phi_x||_\theta^2+(1-\frac{\gamma^2}{4})||\phi||^2_\theta$, we obtain

\begin{equation}\label{ta11}\begin{split}
\frac12 (\frac{dI}{dt}   +  \gamma I) \leq
-{\rm Re}\int \psi_x (1-\theta_R) y+\frac 1 R \int_\R \theta'(\frac x R) y \phi_x+ \int_\R g (1-\theta_R) y.
\end{split}\end{equation}

\noindent We have then

\begin{equation}\label{ta12}
\frac{dI}{dt}   + \frac \gamma 2 I\leq \frac{8}{\gamma}(||\psi_x||^2_{|x|>R}+||g||^2_{|x|>R})+\frac K R.
\end{equation}

\noindent Due to Lemma \ref{tail1} the conclusion follows promptly.

\end{proof}

\subsection{Checking the assumptions of Theorem \ref{chulas}}
\label{secondstepbeyond}

Consider $(\psi,\phi,\phi_t)$ a complete trajectory that belongs to the global attractor.
Introduce $(u,v,v_t)$ that are the differential of $S(t)$ along the trajectory.
Then $(u,v,v_t)$ solves the system, setting $z=v_t+\frac \gamma 2 v$

\begin{eqnarray}
 \label{f1.1} u_t +i\Delta u  +\alpha u  &=&  i\phi u +i \psi v, \\
 \label{f1.2} v_t+\frac \gamma 2 v & = & z, \\
 \label{f1.3} z_t+\frac \gamma 2 z -\Delta v + (1-\frac{\gamma^2}{4}) v  &=&  -Re u_{x}.
\end{eqnarray}

We begin with an $L^2$ estimate. Considering the scalar product of \eqref{f1.1} with $u$ yields,
using $1= (1-\theta_R)+\theta_R$,

\begin{eqnarray}\label{che1}
\frac12 \frac{d}{dt}||u||^2_{L^2}+\alpha ||u||^2_{L^2}&=&-{\rm Im}\int (1-\theta_R) \psi v \overline u- {\rm Im}\int \theta_R \psi v \overline u
\\ &\leq & (||\psi||_\theta ||u||_{H^1(\R)}+K||u||_{L^2(-2R,2R)}) ||v||_{L^2(\R)}\\ &\leq& \frac{\eta}{8}||v||_{L^2(\R)}^2+
\frac 2 \eta||\psi||^2_\theta ||u||_{H^1(\R)}^2+\frac 2 \eta K ||u||_{L^2(-2R,2R)}^2,
\end{eqnarray}

\noindent for a given $\eta>0$ that will be specified in the sequel. Introduce now

$$H(u)=||u_x||^2_{L^2}+\int_\R \phi |u|^2+2 {\rm Re}\int_\R v \psi \overline u.$$

\noindent Considering the scalar product of \eqref{f1.1} by $-i(u_t+\alpha u)$ yields

\begin{equation}\label{che2}
\frac12 \frac{d}{dt}H(u)+\alpha H(u)=\frac12 \int_\R \phi_t |u|^2 +\alpha {\rm Re}\int_\R v \psi \overline u+{\rm Re}\int_\R v_t \psi \overline u
+{\rm Re}\int_\R v \psi_t \overline u.
\end{equation}

\noindent To bound by above the right hand side of \eqref{che2} we use the $1=(1-\theta_R)+\theta_R$ trick.
For the first term we have

\begin{equation}\label{che3}
\int_\R \phi_t |u|^2 \leq ||\phi_t||_\theta ||u||^2_{H^1(\R)} + ||\phi_t||_{L^\infty} ||u||^2_{L^2(-2R,2R)}.
\end{equation}

\noindent We also have, using $u\psi=u\psi(1-\theta_R)+u\psi\theta_R$

\begin{equation}\label{che4}\begin{split}
 {\rm Re}\int_\R v \psi \overline u\leq ||v||_{L^2(\R)}||u\psi||_{L^2(\R)}\leq ||v||_{L^2(\R)}(||\psi||_\theta ||u||_{H^1(\R)}+K||u||_{L^2(-2R,2R)})\\
 \leq \frac \eta 8 ||v||_{L^2(\R)}^2 +2 \frac{||\psi||^2_\theta ||u||_{H^1(\R)}^2}{\eta}+2 \frac K \eta ||u||_{L^2(-2R,2R)}^2,
\end{split}\end{equation}

\noindent and similarly

\begin{equation}\label{che5}
{\rm Re}\int_\R v_t \psi \overline u \leq \frac \eta 8 ||v_t||_{L^2(\R)}^2+2\frac{||\psi||^2_\theta ||u||_{H^1(\R)}^2}{\eta}+2 \frac K \eta ||u||_{L^2(-2R,2R)}^2,
\end{equation}

\noindent and

\begin{equation}\label{che6}
 {\rm Re}\int_R v \psi_t \overline u \leq  \frac \eta 8 ||v||_{L^2(\R)}^2 +2 \frac{||\psi_t||^2_\theta ||u||_{H^1(\R)}^2}{\eta}+2 \frac K \eta ||u||_{L^2(-2R,2R)}^2.
\end{equation}

Fix $\varepsilon>0$ small. Consider $R$ given by Lemma \ref{tail1} and \ref{tail2}.
Introduce now $$H_\kappa(u)=H(u)+\kappa ||u||^2_{L^2},$$
\noindent where $\kappa>0$ will be chosen subsequently depending on $\eta$,
$\varepsilon$ and on the data of the equation. Gathering estimates \eqref{che1}-\eqref{che6} leads to

\begin{equation}\label{dim1}
\frac12 \frac{d}{dt}H_\kappa(u)+\alpha H_\kappa(u)\leq \frac{\varepsilon}{\eta}||u||^2_{H^1}+ \frac \eta 2( ||v||^2_{L^2}+||v_t||^2_{L^2})+ \frac K \eta ||u||^2_{L^2(-2R,2R)}.
\end{equation}

\noindent Introduce now a norm $N$ in the energy space that reads for $V=(\psi,\phi,\phi_t)$,
for $\beta$ a small number that will be specified subsequently (and for $\kappa$ large enough)
\begin{equation}\label{che7}
N(V)^2=\beta ||v_x||^2_{L^2}+(1-\frac{\gamma^2}{4})\beta||v||_{L^2}^2+\beta||z||^2_{L^2}+H_\kappa(u).
\end{equation}

Set
 $$\tilde I= ||z||_{L^2(\R)}^2+||v_x||_{L^2(\R)}^2+(1-\frac{\gamma^2}{4})||v||_{L^2(\R)}^2.$$

\noindent Considering the scalar product of \eqref{f1.2} with $-\Delta v + (1-\frac{\gamma^2}{4})v$ and the scalar
 product of \eqref{f1.3} with $z$ and summing the resulting equations leads to

 \begin{equation}\label{che9}\begin{split}
\frac12\frac{d\tilde I}{dt} + \frac \gamma 2 \tilde I
  \leq -{\rm Re}\int_\R u_x z
\leq \frac{\gamma}{4} ||z||^2_{L^2}+ \frac{1}{\gamma} ||u_x||_{L^2}^2.
\end{split}\end{equation}

\noindent Multiplying this inequality by $\beta$ and adding \eqref{dim1} we have

\begin{equation}\label{dim2}\begin{split}
\frac 12\frac{d}{dt} N(V)^2 + \min(\alpha, \frac \gamma 4)N(V)^2\leq \\
(\frac \varepsilon \eta+\frac \beta \gamma)||u||^2_{H^1}+\frac \eta 2(||v||^2_{L^2}+||v_t||^2_{L^2})+\frac K \eta ||u||^2_{L^2(-2R,2R)}.
\end{split}\end{equation}

We now chose the parameters as follows. Let us observe that by Cauchy-Schwarz's inequality

$$ H_\kappa(u)\geq ||u_x||^2_{L^2}-\frac \beta 4 ||v||^2_{L^2}+(\kappa -\frac {\tilde K}{\beta} )||u||^2_{L^2}.$$

\noindent We set $\kappa=1+\frac {\tilde K}{\beta}$. Therefore $N(V)^2\geq ||u||^2_{H^1}+\frac \beta 4 \tilde I$
is the square of a norm. We now fix $\beta$ small enough such that $ 10 \frac \beta \gamma \leq \min(\alpha, \frac \gamma 4),$
then $\eta$ small enough such that $10 \eta \leq \beta \min(\alpha, \frac \gamma 4)$
and finally $\varepsilon$ small enough such that $10 \frac{K\varepsilon}\eta \leq \min(\alpha, \frac \gamma 4)$.

\noindent We then infer from \eqref{dim2}

\begin{equation}\label{dim3}
\frac{d}{dt} N(V)^2 + \min(\alpha, \frac \gamma 4)N(V)^2\leq \frac K \eta ||u||^2_{L^2(-2R,2R)}.
\end{equation}

\noindent This completes the proof of \eqref{boa}.

\section{Miscellaneous remarks on the global attractor}

What can we say about the structure of the global attractor ?

\subsection{Trivial dynamics}
If $\mathcal{A}$ is reduced to one point, i.e. one stationary solution,
and all solutions converge to this stationary solution, the dynamics are trivial. This is the case if $f=0$
or if $\alpha$ and $\gamma$ are large enough; more precisely

\begin{proposition}
Consider a stationary solution $(\psi^*,\phi^*,0)$.
If $\alpha\gamma>2||\psi^*||_{L^\infty}$ then all solutions converge to
the stationay solution when $t$ diverges to $+\infty$.
\end{proposition}

\begin{proof}

Set $u=\psi-\psi^*$ and $v=\phi-\phi^*$. Then

\begin{eqnarray}
\label{pog.1} u_t+\alpha u+ iu_{xx} & = & i\phi^* u +i v \psi^*, \\
\label{pog.2} v_{tt}+\gamma v_t +v -\Delta v & = & -{\rm Re}u_x.
\end{eqnarray}

\noindent Since

$$ |\int_\R v \psi^*\overline u|\leq ||\psi^*||_{L^\infty} ||v||_{L^2(\R)} ||u||_{L^2(\R)},$$

\noindent considering the scalar product of \eqref{pog.1} with $u$ leads to

\begin{equation}\label{pog.3}
\frac{d}{dt} ||u||_{L^2(\R)}+\alpha ||u||_{L^2(\R)} \leq ||\psi^*||_{L^\infty} ||v||_{L^2(\R)}.
\end{equation}

\noindent Using the representation formula \eqref{rep} yields

\begin{equation}\label{pog.4}
||v(t)||_{L^2(\R)}\leq \int_{-\infty}^t e^{-\frac \gamma 2(t-s)} ||u(s)||_{L^2(\R)}ds.
\end{equation}

\noindent Introducing $q=\int_{-\infty}^t e^{-\frac \gamma 2(t-s)} ||u(s)||_{L^2(\R)}ds$
we then infer from \eqref{pog.3}, \eqref{pog.4} that

\begin{eqnarray}
\label{pog.5} \frac{d}{dt} ||u||_{L^2(\R)}+\alpha ||u||_{L^2(\R)}& \leq & q||\psi^*||_{L^\infty} \\
\label{pog.6} \dot q +\frac \gamma 2 q &= &||u||_{L^2(\R)}.
\end{eqnarray}

\noindent The condition $\alpha\gamma>2||\psi^*||_{L^\infty}$ ensures that $u$, $q$
and also $v$ goes to $0$ when $t$ diverges to $+\infty$.
\end{proof}

\subsection{An example of non trivial dynamics}

We provide below a concrete example of KGS system with at least
two stationary solutions. We proceed as follows.
Assume $g=A_x$ is the derivative of a smooth function $A$ in $H^1(\R)$.
We seek first a solution of \eqref{e.1}, \eqref{e.2} such that $\phi^*=0$.
This amounts to solve, setting $a^*=$Re$\psi^*$ and $b^*=$Im$\psi^*$,

\begin{eqnarray}
a^*&=&A, \\
\label{s.1} -\Delta b^*  +\alpha A &= & {\rm Re} f, \\
\label{s.2} \Delta A + \alpha b^* & =& {\rm Im} f.
\end{eqnarray}

\noindent This leads to the compatibility condition involving $A$ and the
real and imaginary parts of $f$

\begin{equation}\label{s.3}
\Delta^2 A + \alpha^2 A= \alpha {\rm Re} f+ \Delta {\rm Im} f.
\end{equation}

We seek another stationary solution that satisfies $a^*=0$.
Introducing $k(x)=\frac12 \exp(-|x|)$, equation \eqref{e.2} yields

\begin{equation}\label{s.4}
\phi^*(x)= k*A_x(x)= \frac12 (e^x\int_x^{+\infty} e^{-y} A(y)dy- e^{-x} \int_{-\infty}^x e^y A(y)dy).
\end{equation}

\noindent Equation \eqref{e.1} leads to

\begin{eqnarray}
\label{s.5} -\Delta \tilde b^*+\phi^* \tilde b^* &= & {\rm Re} f, \\
\label{s.6}  \alpha \tilde b^* & =& {\rm Im} f.
\end{eqnarray}

\noindent This leads to another compatibility condition that reads

\begin{equation}\label{s.7}
-\Delta {\rm Im} f + \phi^* {\rm Im f} =\alpha {\rm Re} f.
\end{equation}

\noindent Comparing \eqref{s.4} and \eqref{s.7} we can find two different solutions, for a given $\alpha$,
if the pair $(A, {\rm Im} f)$ solves the integro-differential equation

\begin{equation}\label{s.8}
\Delta^2 A + \alpha^2 A= {\rm Im} f (k*A)_x,
\end{equation}

\noindent and $A\neq 0$. We now choose $A(x)=\exp(-\frac{x^4}{4})$. Roughly speaking
we aim to set

\begin{equation}\label{rough}
{\rm Im} f= \frac{\Delta^2 A + \alpha^2 A}{(k*A)_x}.
\end{equation}

This is possible if

\begin{itemize}
\item The zeros of $(k*A)_x$ are also zeros of $\Delta^2 A + \alpha^2 A$ with at most the same order.
\item The decay at $\pm \infty$ of $(k*A)_x$ is much slower than the decay of $\Delta^2 A + \alpha^2 A$, in order to ensure that Im$f$ is in $L^2(\R)$.
\end{itemize}

We now state

\begin{lemma}\label{fin}
The function $(k*A)_x$ has only one zero in $\R$, that is $0$ and that is simple.
The function $(k*A)_x$ decays as O$(\exp-|x|)$ at the infinity.
\end{lemma}

Postponing the proof of Lemma \ref{fin} we observe that
$$ \Delta^2 A(0) + \alpha^2 A(0)= -6 +\alpha^2. $$
We now fix $\alpha=\sqrt 6$ in order $\Delta^2 A + \alpha^2 A$ vanishes at $0$.
Since $\Delta^2 A + 6 A \simeq cx^4 \exp(-\frac{x^4}{4})$ then \eqref{rough} is valid and we have two different
stationary solutions for this $\alpha, A, f$. We now prove Lemma \ref{fin}.

\begin{proof}
The function $(k*A)_x$ is odd, then it vanishes at $0$.
Introducing $W=k*A$ we have that $W(0)<A(0)=1$ and since $-\ddot W+W=A$ then $V'(0)\neq 0$.
Moreover the representation formula, for $x>0$,

\begin{equation}\label{s.10}
2 (k*A)_x(x)=\int_0^{+\infty} (e^{-|x-y|}-e^{-x-y})A'(y)dy,
\end{equation}

\noindent ensures that the only zero of $(k*A)_x$ is $x$.

Besides, we easily infer from the representation formula \eqref{s.4}
that $$(k*A)_x(x) \simeq -\frac{e^{-x}}{2}(\int_\R e^{y} A(y)dy),$$
\noindent when $x$ goes to $+\infty.$

\end{proof}

\subsection{Conclusion}

The structure of the global attractor is known when there exists a Lyapunov
function. In this case the global attractor is the set of stationary solutions
and the heteroclinic orbits that connects the stationary solutions.
Some dissipative perturbation of Hamiltonian systems have a Lyapunov function.
For parametric KGS equation this is not the case; when reduced to the attractor,
the dynamics is the dynamics of a Schr\"odinger equation with an infinite memory term.
This point of view will be addressed in a forthcoming work.

\section*{Acknowledgments} The second author acknowledges the hospitality of the Universit\'e
de Lille. Both authors acknowledge the support of the CDP C2EMPI, together with the French State under the France-2030 programme, the University of Lille, the Initiative of Excellence of the University of Lille, the European Metropolis of Lille for their funding and support of the R-CDP-24-004-C2EMPI project

\medskip
Received XXXXXX; revised XXXXXX.
\medskip

\end{document}